\documentclass{amsart}
\usepackage[cp1251]{inputenc}
\usepackage[english]{babel}

\usepackage{amsmath}
\usepackage{amssymb}
\usepackage{amsfonts}

\begin{document}

\title[Explicit solution of the problem of equivalence \ dots ]
{Explicit solution of the problem of equivalence for some Painleve equations}

\author{ V.V.Kartak }\footnote{\uppercase{W}ork partially
supported by grant  08-01-97020-p\_povoljie\_a of the \uppercase{R}ussian 
\uppercase{F}ond of \uppercase{B}asic \uppercase{R}esearch, coordinator Prof. Fazullin Z.Yu.}

\address{Bashkir State University \\
Department of Mathematics\\
450074 Frunze Street\\ 
Ufa, Russia\\ 
E-mail: kvera@mail.ru}


\maketitle 
{
\small
\begin{quote}
{\bf Abstact.} 
For an arbitrary ordinary second order  differential equation
 a test is constructed that checks if this equation is equivalent to
Painleve I, II or Painleve III with three zero parameters equations
under the substitutions of variables. If it is true then in case the
Painleve equations I and II an explicite change of variables is given
that is written using the differential invariants of the equation.

\medskip
 
 \noindent{\bf Keywords:}  Painleve equations, equivalence problem, differential invariants.
 \end{quote}
}

\section{Introduction}

It is well-known fact that the ordinary differential equations of the form 
\begin{equation}\label{eq}
y''=P(x,y)+3\,Q(x,y)y'+3\,R(x,y)y^{\prime 2}+S(x,y)y^{\prime 3}
\end{equation}
preserve their form
 under the action of arbitrary point transformations 
\begin{equation}\label{zam}
\tilde x=\tilde x(x,y),\quad \tilde y=\tilde y(x,y).
\end{equation}
So we can use of geometric methods to study these equations.
See~ \cite {Liouville}, \cite{Lie}, \cite{Tresse1}, \cite{Tresse2}, \cite{Cartan}.
In particular, we can  build differential invariants associated with
the equation (\ref{eq}) from its coefficients - functions $ P (x, y)
$, $ Q (x, y) $,
$ R (x, y) $, $ S (x, y) $, and their derivatives.
Such invariants are called {\it Cartan invariants}. 
See \cite{Thomsen}, \cite{Grissom}, \cite{Rom}, \cite{BordagDruima}.

All six famous Painleve equations (\cite{Its}, \cite{Sigur}, \cite{Gromak}) have the form (\ref{eq}). The first, second and third Painleve equations respectively are:
$$
 y''=6 y^2+ x,\qquad   y''=2 y^3+ x y+a,\qquad 
y^{\prime\prime}=\frac     1y(y^{\prime})^2-\frac     1x
y^{\prime}+\frac 1x(ay^2+b)+cy^3+\frac dy.
$$
We use methods of differential invariants, described in the papers \cite{Sharipov1}, \cite{Sharipov2}, \cite{Sharipov3} in order to solve the problem of the equivalence for Painleve I and II equations and for Painleve III equation with three zero parameters. 

Similar studies were conducted previously, see \cite{Lamb},
\cite{Sokolov}, \cite{Babich},  \cite{Hietarinta}, \cite{Dridi}.
 However, for the first time the check test for the equivalence is
formulated in the effectively verifiable manner, it can be programmed. This is a continuation of work \cite{Kartak}.

For further calculations we also need to define pseudotensorial field and its covariant derivative (as they were formulated in \cite{Sharipov1}). 

{\bf Definition 1. }{ \it The pseudotensorial field weights $ m $ of valence $ (r, s) $ is called indexed set of variables that  transform with the following rule under change of coordinate system}
$$
F^{i_1\dots i_r}_{j_1\dots j_s}=(\det T)^m{\sum_{p_1\dots p_r}}
{\sum_{q_1\dots q_s}} S^{i_1}_{p_1}\dots S^{i_r}_{p_r}T^{q_1}_{j_1}
\dots T^{q_s}_{j_s} \tilde F^{p_1\dots p_r}_{q_1\dots q_s},
$$
here $ T $ is the inverse transfer matrix  when one changes one
coordinate system to another system in the plane.
We see that there is the only factor $ (\ det T)^ m $  that
distinguishes  Definition 1 from the classical definition of tensorial
field.

{\bf Definition 2.}{\it  
The following object is called the covariant derivative of
pseudotensorial field $ F $ of valence $ (r, s) $ and weight $ m $:
}
$$
\aligned
\nabla_kF^{i_1\dots i_r}_{j_1\dots j_s}=& \frac{\partial F^{i_1\dots i_r}_{j_1\dots j_s}}{\partial x^k}+
\sum_{n=1}^r\sum_{v_n=1}^2\Gamma_{k v_n}^{i_n}F^{i_1\dots v_n\dots i_r}_{j_1\dots j_s}-\\
-&\sum_{n=1}^s\sum_{w_n=1}^2\Gamma_{k j_n}^{w_n}F^{i_1\dots i_r}_{j_1\dots w_n\dots  j_s}+
m\varphi_k F^{i_1\dots i_r}_{j_1\dots j_s}
\endaligned
$$
Under the covariant differentiation the pseudotensorial field $ F $ of valence $ (r, s) $ and weight $ m $ becomes the pseudotensorial field of valence $ (r, s +1) $ and weight $ m. $ 

Only the term $ m \varphi_k F ^ (i_1 \dots i_r) _ (j_1 \dots j_s) $
distinguishes  Definition 2
from the definition of covariant derivative of tensor fields.
Here $ \varphi_1 $ and $ \varphi_2 $ are auxiliary coefficients,
the explicit formula (\ref{phi1}), (\ref{phi2}) for their calculation
are contained in the Appendix.

\section{Painleve I equation }

Suppose we have a certain equation  (\ref{eq}). We are looking for a change of variables, which translates it to the Painleve I equation
\begin{equation}\label{pen1} 
\tilde y''=6\tilde y^2+\tilde x.
\end{equation}

There are two  pseudovectorial fields $ \alpha $ of weight 2 and $ \theta $ of weight -1 
 associated with equation (\ref{pen1})
(for details, see \cite{Sharipov2}). Explicit formulas (\ref{alpha}), (\ref{theta}) for calculating the
coordinates of these fields are contained in the Appendix.
 
For the equation (\ref{pen1}) they are given by: 
$$
\tilde \alpha^1=\tilde B=0,\quad \tilde\alpha^2=- \tilde A=-12,\quad \tilde\theta^1=-\frac{1}{12},\quad \tilde\theta^2=0.
$$
Then the transformation laws of these fields under the change of
coordinates ( \ref{zam}) are the following:
\begin{equation}\label{p1}
\begin{aligned}
 & \left(\aligned B\\-A\endaligned\right)=\det T \left(\begin{array}{cc}
	 \tilde y_{0.1} & -\tilde x_{0.1}\\ -\tilde y_{1.0} & \tilde x_{1.0}
	 \end{array}\right)  \left(\aligned 0\\12\endaligned\right),\\
 & \left(\aligned \theta^1\\ \theta^2\endaligned\right)=
\frac 1{\det^2 T} \left(\begin{array}{cc}
	 \tilde y_{0.1} & -\tilde x_{0.1}\\ -\tilde y_{1.0} & \tilde x_{1.0}
	 \end{array}\right)  \left(\aligned  -{1}/{12}\\ 0 \endaligned\right).
	 \end{aligned}
\end{equation}

In addition, the coefficients of the equation (\ref{pen1}) define the pseudoinvariants, explicit formulas for the calculation of which are contained in the Appendix. 
$ N $ of weight 2 (\ref{N}); $ \Omega $ of weight 1 (\ref{Omega1}), ( \ref{Omega2}); $ \Theta $ of weight -2 (\ref{Theta}); $ L $ of weight -4 (\ref{L}); 
$ L_1 $ weight -5 (\ref{L1}); $ V $ of weight -3 (\ref{V}); $ W $ of weight -6 (\ref{W}). 
Now we can construct the set of invariants: 
$$
I_1=\frac{L_1^4}{L^5},\qquad  I_2=\frac{\Theta^2}{L}.
$$
In this case pseudoinvarians are:
$$
\tilde N=0,\quad \tilde \Omega=0,\quad \tilde \Theta =-\frac{\tilde y}{12},\quad \tilde L=\frac{\tilde x}{12^3},\quad \tilde L_1=-\frac{1}{12^4},\quad
\tilde V=0,\quad \tilde W=0.
$$

We see that pseudoinvariant $ \tilde L_1 $ of weight -5 for the
Painleve equation I (\ref{pen1}) is a constant.

We can calculate the  determinant of the inverse matrix of transition
$ \det T $ using the transformation law for the $ \tilde L_1 $ under
the change of coordinates:
$$
L_1=\frac{ \tilde L_1}{(\det T)^5}=\frac{-1}{12^4(\det T)^5},\quad \det T= \sqrt[5]{-\frac{1}{12^4L_1}}.
$$
Now from (\ref{p1}) we calculate partial derivatives of the coordinate
functions:
$$
\tilde x_{1.0}=-A\sqrt[5]{-\frac{L_1}{12}},\quad
\tilde x_{0.1}=-B\sqrt[5]{-\frac{L_1}{12}},\quad
\tilde y_{1.0}=\frac{\theta^2}{\sqrt[5]{12^3L_1^2}},\quad
\tilde y_{0.1}=-\frac{\theta^1}{\sqrt[5]{12^3L_1^2}}.
$$
For any equation (\ref{eq}) functions $ A $, $ B $, $ \theta ^ 1 $, $
\theta ^ 2 $, $ L_1 $ are known, they are determined by functions $ P,
$ $ Q $, $ R, $ $ S $ via explicite formulae (\ref{alpha}),
(\ref{theta}), (\ref{L1}) given in the Appendix.
From the last formula it is easy to obtain conditions for the compatibility:
$$
\left(A\sqrt[5]{-\frac{L_1}{12}} \right)_y=\left( B\sqrt[5]{-\frac{L_1}{12}}  \right)_x,\quad
\left(\frac{\theta^2}{\sqrt[5]{12^3L_1^2}} \right)_y=\left( -\frac{\theta^1}{\sqrt[5]{12^3L_1^2}}\right)_x.
$$

The first condition of compatibility gives us the following expression 
$$
\aligned
&5L_1(A_{y}-B_{x})+A(L_1)_y-B(L_1)_x=5L_1(A_{y}-B_{x})-\nabla_{\alpha}L_1+5L_1(\varphi_1B-\varphi_2A)=\\
&=5L_1(A_{y}-B_{x}-\varphi_1B+\varphi_2A)-V=
6L_1 N-V=0.
\endaligned
$$
We have used the definition of pseudoinvariant $V$
$$
V=\nabla_{\alpha}L_1= (L_1)_{x}B- (L_1)_{y}A-5L_1(B\varphi_1-A\varphi_2),
$$
as well as the following clearly verifiable identity: 
\begin{equation}\label{1}
B_x-A_y=-\frac 65 N+\varphi_2A-\varphi_1B.
\end{equation}

From the second condition of compatibility:
$$
\aligned
&\quad\frac 52 L_1\left(\theta^1_x+\theta^2_y\right)-\left((L_1)_x\theta^1+(L_1)_y\theta^2\right)
=\\
&=\frac 52 L_1\left(\theta^1_x+\theta^2_y\right)-\nabla_{\theta}L_1-5L_1(\varphi_1\theta^1+\varphi_2\theta^2)=\\
&=\frac 52 L_1\left( (\Theta_{ y}-2\varphi_2\Theta)_x+(-\Theta_{
x}+2\varphi_1\Theta)_y      \right)- W-\\
&-5L_1\left(\varphi_1(\Theta_{y}-2\varphi_2\Theta)+
\varphi_2(-\Theta_{x}+2\varphi_1\Theta)\right)=\\
&=\frac 52 L_1\left(\left(\Theta_{xy}-2(\varphi_2)_x\Theta-2\varphi_2\Theta_{x}\right)+\left(-\Theta_{xy}+
2(\varphi_1)_y\Theta+2\varphi_1\Theta_{y}\right)\right)- W-\\
&-5L_1\left(\varphi_1\Theta_{y}-
\varphi_2\Theta_{x})\right)
=5L_1((\varphi_1)_y-(\varphi_2)_x)\Theta-W=3L_1\Omega\Theta-W=0.
\endaligned
$$
We have used the definition of pseudoinvariants $W$ and $\Omega$, pseudovectorial field $\theta$:
$$
W=\nabla_{\theta}L_1= (L_1)_{x}\theta^1+ (L_1)_{
y}\theta^2-5L_1(\varphi_1\theta^1+\varphi_2\theta^2),\quad \Omega=\frac 53 \left((\varphi_1)_y-(\varphi_2)_x\right),
$$
$$
\theta^1=\Theta_{ y}-2\varphi_2\Theta, \qquad \theta^2=-\Theta_{
x}+2\varphi_1\Theta.
$$

As for the Painleve I equation all pseudoinvariants $ \tilde N $, $ \tilde \Omega $, $ \tilde V $, 
$  \tilde W $ are identically equal to zero, 
they must be zero for any equation that is equivalent to Painleve I equation.

Since $ N = 0 $ and $ V = 0 $, the first condition of compatibility is true,
and since $ W = 0 $ and $ \Omega = 0 $, the second condition of
compatibility is true.
Thus, realisation of these conditions is equal to the existing of the point substitution of variables, which reduces original equation (\ref{eq}) to equation (\ref{pen1}).

Invariants of the equation ( \ref{pen1}) are given by:
$$
I_1=\frac{1}{12\tilde x^5},
\quad I_2=\frac{12 \tilde y^2}{\tilde x}.
$$
Let us resolve $I_1$ and $I_2$ relatively $\tilde x$ and $\tilde y$.
The explicite change of variables is given by (\ref{p1zam})
 and find the explicite change of variables:\begin{equation}\label{p1zam}
\tilde x=\frac 1{\sqrt[5]{12 I_1}}, \qquad \tilde y=\pm\frac{\sqrt{I_2}}{\sqrt[5]{12^3}\sqrt[10]{I_1}}.
\end{equation}

{\bf Theorem 1} {\it Equation (\ref{eq}) is equaivalent to Painleve I equation (\ref{pen1}) under trans\-for\-mations (\ref{zam}) if and only if the following conditions are true:
1) $F=0$ (\ref{F}), but  $A\ne 0$ or $B\ne 0$ (\ref{alpha}),
2) $\Omega=0$ (\ref{Omega1}), (\ref{Omega2}),
3) $N=0$ (\ref{N}), 
4) $W=0$ (\ref{W}),
5) $V=0$ (\ref{V}),
6) $\Theta\ne 0$ (\ref{Theta}),
7) $L_1\ne 0$ (\ref{L1}).
The explicite change of variables is given by (\ref{p1zam}).}

{\bf Example 1.}  The following equation is equivalent to the Painleve I equation:
$$
\aligned
y''=&-\sin^3 y(6x\cos^2 y+\sin y)+
\frac 1x (-18x^3\cos^3 y\sin^2 y-3x^2\sin^3 y\cos y-2)y^\prime-\\
&-(18x^3\cos^4 y\sin y+3x^2\sin^2 y\cos^2 y)y^{\prime 2}
 -(6x^4\cos^5 y+x^3\sin y\cos ^3y+x)y^{\prime 3}.
 \endaligned
 $$
 Invariants and the explicite change of variables are the following:
 $$
I_1=\frac 1{12}\frac{1}{x^5\sin^5 y}, \quad I_2=\frac{12x\cos^2 y}{\sin y},
\quad \tilde y=x\cos y, \quad \tilde x=x\sin y.
$$

\section{Painleve II equation}

For a certain equation of type (\ref{eq}) we are looking for a change
of variables that transformes it into the Painleve II equation
\begin{equation}\label{pen2} 
\tilde y''=2\tilde y^3+\tilde x\tilde y+a.
\end{equation}

There are two pseudovectorial fields  $\alpha$ of weight  2
(\ref{alpha}) and $\xi$ of weight  3 (\ref{xi}) associated with
equation (\ref{pen2}).
 Explicit formulas for their calculation are contained in Appendix. For equation (\ref{pen2}) they are:
$$
\tilde \alpha^1=\tilde B=0,\quad \tilde\alpha^2=- \tilde A=-12\tilde y,\quad \tilde\xi^1=-\frac{24}{5\tilde y},\quad \tilde\xi^2=0.
$$
Under the change of variables they transform into the rule:
$$
 \left(\aligned B\\-A\endaligned\right)=(\det T) \left(\begin{array}{cc}
	 \tilde y_{0.1} & -\tilde x_{0.1}\\ -\tilde y_{1.0} & \tilde x_{1.0}
	 \end{array}\right)  \left(\aligned  0\\-12\tilde y\endaligned\right),
$$
$$
 \left(\aligned \xi^1\\ \xi^2\endaligned\right)=(\det T)^2 \left(\begin{array}{cc}
	 \tilde y_{0.1} & -\tilde x_{0.1}\\ -\tilde y_{1.0} & \tilde x_{1.0}
	 \end{array}\right)  \left(\aligned  -{24}/{(5\tilde y)}\\ 0 \endaligned\right).
 $$
 
Pseudoinvariants $M$ of weight 4 (\ref{M1}), (\ref{M2}), $N$ of weight 2 (\ref{N}), $\Omega$ of weight 1 (\ref{Omega1}), (\ref{Omega2}) 
and $\Gamma$ (\ref{Gamma}) for the equation (\ref{pen2}) are given by:
  $$\tilde M=\frac{288}{5},\qquad\tilde N=4,\qquad \tilde\Omega=0,\qquad \tilde\Gamma=\frac {48}{25}\frac{2\tilde y^3+\tilde x\tilde y+a}{\tilde y^3}.$$

Invariants of the equation are calculated by the formulas:
\begin{equation}\label{i1}
\aligned
& I_1=\frac M{N^2},\quad I_3=\frac {\Gamma}{M},\quad I_6=\frac{\nabla_{\alpha}I_3}{N}=\frac{B(I_3)'_x-A(I_3)'_y}{N},\\
& I_9=\frac{(\nabla_{\gamma}I_3)^2}{N^3}=\frac{(\xi^1(I_3)'_x+\xi^2(I_3)'_y)^2}{N^3}.
\endaligned
\end{equation}
Similarly to the previously discussed case of the conversion formula for $ N $
let us find $ \det T: $ 
$$
N=4(\det T)^2,\qquad \det T=\frac{\sqrt N}{2},
$$
then
\begin{equation}\label{2}
\frac{\tilde y_{0.1}}{\tilde y}=-\frac{5}{6}\frac{\xi^1}{N},\quad
\frac{\tilde y_{1.0}}{\tilde y}=\frac{5}{6}\frac{\xi^2}{N} 
\end{equation}
and the corresponding compatibility condition has the form:
$$
\left(-\frac{\xi^1}{N}\right)_x=\left(\frac{\xi^2}{N}\right)_y.
$$
It is equivalent to
$$\aligned
& N(\xi^1_x+\xi^2_y)-(\xi^1 N_x+\xi^2 N_y)=N(N_y+2\varphi_2N)_x+N(-N_x-2\varphi_1N)_y-\\
&-N_x(N_y+2\varphi_2N)-N_y(-N_x-2\varphi_1N)=2N^2((\varphi_2)_x-(\varphi_1)_y)=
\frac {10}3N^2\Omega=0.
\endaligned
$$
We have used the definition
\begin{equation}\label{3}
\xi^1=N_y+2\varphi_2N,\qquad \xi^2=-N_x-2\varphi_1N.
\end{equation}
This condition is fulfilled if $\Omega=0$
As we prove the existence of $ \tilde y $, let's substitute it into the first equality: 
$$
\tilde x_{0.1}=\frac{B}{6\tilde y\sqrt{N}},\quad
\tilde x_{1.0}=\frac{A}{6\tilde y\sqrt{N}},\quad 
\left(\frac{B}{\tilde y\sqrt{N}}\right)_x=\left(\frac{A}{\tilde y\sqrt{N}}\right)_y.
$$
Let us write this in more details:
$$
\aligned
& \frac{B_x-A_y}{\tilde y \sqrt{N}}+\frac{A\tilde y_{0.1}-B\tilde y_{1.0}}{\tilde y^2 \sqrt{N}}+\frac{AN_y-BN_x}{2\tilde y N\sqrt{N}}=\\
&=\frac{-\frac 65N+\varphi_2 A-\varphi_1 B}{\tilde y\sqrt{N}}-\frac{5(B\xi^2+A\xi^1)}{6\tilde y N\sqrt{N}}+
\frac{A(\xi^1-2\varphi_2N)+B(\xi^2+2\varphi_1N)}{2\tilde y N\sqrt{N}}=\\
&=\frac{-6N}{5\tilde y\sqrt{N}}-\frac{B\xi^2+A\xi^1}{3\tilde y N\sqrt{N}}=
\frac{-18 N^2+5 M}{15\tilde y N\sqrt{N}}=\frac{-18+5I_1}{15\tilde y N^3\sqrt{N}}=0.
\endaligned
$$
The first condition of compatibility is satisfied if
  $I_1=18/5$.
We used formula (\ref{1}), (\ref{2}), (\ref{3}) and the definition
$$
M=-A\xi^1-B\xi^2.
$$

Values of basic invariants (\ref{i1}) to the equation (\ref{pen2}): 
 $$
\aligned
 I_1&=\frac{18}{5},\quad I_3=\frac{2\tilde y^3+\tilde x\tilde y+a}{30\tilde y^3},\quad I_6=\frac{2\tilde
 x\tilde y+3a}{10\tilde y^3},\quad I_9=\frac{1}{2500\tilde y^6}.
\endaligned
$$
Let us construct a  new invariant, which is up to the sign equals the
parameter of equation (\ref{pen2}):
\begin{equation}\label{j}
J=\frac{1}{50}\frac{4+10I_6-60I_3}{\sqrt{I_9}}=\pm a.
\end{equation}

{\bf Lemma 1.} {\it Equations Painleve II with the different parameters $a_1\ne\pm a_2$
are non-equivalent.}

Via the formula of the invarians we find the explicite change of variables:
\begin{equation}\label{p2zam}
\tilde y=\frac{1}{\sqrt[6]{2500 I_9}},\qquad
\tilde x=\frac{5I_6}{\sqrt[6]{2500I_9}}-\frac 32 J\sqrt[6]{2500 I_9}.
\end{equation}

{\bf Theorem 2.} 
{\it An arbitrary equation (\ref{eq}) is equivalent to the Painleve II equation with parameter $ a = \pm J $ (\ref{j}) if and only if the following conditions are true:
1) $F=0$ (\ref{F}), but $A\ne 0$ or $B\ne 0$ (\ref{alpha}),
2) $\Omega=0$ (\ref{Omega1}), (\ref{Omega2}) ,
3) $M\ne 0$ (\ref{M1}), (\ref{M2}), 4) $I_1=18/5$ (\ref{i1}). The explicit change of variables is (\ref{p2zam}).}

{\bf Example 2.} Under a linear change of variables equation 6.9 from \cite{Kamke}
 is reduced to the Painleve II equation with parameter $ \pm J $:
$$
\aligned
y^{\prime\prime}&=-ay^3-bxy-cy-d.\\
J&=-\sqrt{\frac a2}\cdot\frac{d}{b},\qquad
\tilde y=\sqrt{\frac a2}\cdot\frac{y}{\sqrt[3]{b}}, \qquad
\tilde x=-\frac{bx+c}{\sqrt[3]{b^2}}.
\endaligned
$$

\section{Equation Painleve III with three zero parametrs }

A general form of the Painleve equations III is the following:
$$
y^{\prime\prime}=\frac     1y(y^{\prime})^2-\frac     1x
y^{\prime}+\frac 1x(ay^2+b)+cy^3+\frac dy.
$$
It is a 4-parameter family of equations, which we denote by $PIII (a, b, c, d)$.

If three out of four of these parameters are zero, then these equations of Painleve III have special properties:

1. They have a two-dimensional algebra of point symmetries and hence integrable. See \cite{Gromak}.

2. All these equations are equivalent to each other. Referring to  work \cite{Hietarinta}, we write the change of variables:
$$
\aligned
& PIII(0,b,0,0)\stackrel{1)}{\rightarrow}PIII(-b,0,0,0)\stackrel{2)}{\rightarrow}PIII(0,0,-b,0)\stackrel{3)}{\rightarrow}PIII(0,0,0,b),\\
&\text{here\;\;}1), 3)\:\;x=\tilde x,\,y=1/{\tilde y},\quad 2)\; x=\tilde x^2/2,\,y=\tilde y^2.
\endaligned
$$

Therefore it makes sense to solve the problem of equivalence for the same type of equations. We chose $PIII (0, b, 0,0)$: 
\begin{equation}\label{pen3b}
y^{\prime\prime}=\frac     1y(y^{\prime})^2-\frac     1x
y^{\prime}+\frac bx.
\end{equation}

For the equation (\ref{pen3b}) the coordinates of the pseudovectorial fields $ \alpha $ of weight 2 (\ref{alpha}) and $ \xi $ of the weights 3 (\ref{xi}) are: 
$$
\tilde A=\frac{b}{\tilde x\tilde y^3}, 
\qquad \tilde B=0,\qquad
\tilde\xi^1=-\frac{1}{15}\frac{b}{\tilde x\tilde y^4}, \qquad
\tilde\xi^2=-\frac{1}{15}\frac{b}{\tilde x^2\tilde y^3},
$$
and  values of pseudoinvarians $ M $ of the weight 4 (\ref{M1}), (\ref{M2}), $ N $ of weight 2 (\ref{N}) and $ \Omega $ of weight 1 (\ref{Omega1}) , (\ref{Omega2}) are: 
$$\qquad \tilde N=-\frac 13\frac{b}{\tilde x\tilde y^3},\qquad
\tilde M=\frac{1}{15}\frac{b^2}{\tilde x^2\tilde y^6}, \qquad \tilde\Omega=0.
$$

The basic invariants of the equation are:
$$
I_1=\frac{M}{N^2}=\frac 35,\qquad I_2=\frac{\Omega^2}N=0,\qquad I_3=\frac{\Gamma}{M}=\frac 1{15}.
$$

Let us calculate 
$ \det T $ from the transformation law of the pseudoinvariant $ N $: 
$$
\det T=\frac{\sqrt{-3N\tilde x\tilde y^3}}{\sqrt b},
$$
then from the laws of transformation of the pseudovectorials fields $ \alpha $ and $ \xi $ we get 
\begin{equation}\label{x}
\tilde x_{0.1}=\frac{B\sqrt{\tilde x\tilde y}}{\sqrt{-3bN}}, \qquad
\tilde x_{1.0}=\frac{A\sqrt{\tilde x\tilde y}}{\sqrt{-3bN}},
\end{equation}
\begin{equation}\label{y}
\tilde y_{0.1}=\frac{5\xi^1\tilde y}{N}+\frac{B\tilde y^2}{\sqrt{-3bN\tilde x\tilde y}},\qquad
\tilde y_{1.0}=-\frac{5\xi^2\tilde y}{N}+\frac{A\tilde y^2}{\sqrt{-3bN\tilde x\tilde y}}.
\end{equation}
In this case the compatibility conditions are the following: $$
\left(\frac{B\sqrt{\tilde y}}{\sqrt{N}}\right)_x=\left( \frac{A\sqrt{\tilde y}}{\sqrt{N}} \right)_y,
\quad
\left( \frac{5\xi^1}{N}+\frac{B\sqrt{\tilde y}}{\sqrt{-3bN\tilde x}} \right)_x=
\left( -\frac{5\xi^2}{N}+\frac{A\sqrt{\tilde y}}{\sqrt{-3bN\tilde x}} \right)_y
$$

Let us write the first condition of compatibility, using formulas
(\ref{x}), (\ref{y}):
$$
\aligned
&\frac{(B_{1.0}-A_{0.1})\sqrt{\tilde y}}{\sqrt N}+
\frac{(B\tilde y_{1.0}-A\tilde y_{0.1})}{2\sqrt N\sqrt{\tilde y}}-
\frac 12\frac{\sqrt{\tilde y}(BN_{1.0}-AN_{0.1})}{\sqrt{N^3}}=\\
&=-\frac 65 \sqrt N\sqrt{\tilde y}-\frac 12\frac{\sqrt{\tilde y}M}{\sqrt{N^3}}-
\frac{5\sqrt{\tilde y}}{2\sqrt N}(B\xi^2+A\xi^1)=\sqrt{\tilde y}\sqrt N
\left(-\frac 65-\frac 12 \frac{M}{N^2}+\frac 52 \frac{M}{N^2} \right)=\\
&=\sqrt{\tilde y}\sqrt N
\left(-\frac 65+2 I_1 \right)=0.
\endaligned
$$
It is true if $I_1=3/5.$ The second condition can be written as follows:
$$
\aligned
&\frac{5(\xi^1_x+\xi^2_y)}{N}-\frac{5(\xi^1N_{1.0}+\xi^2N_{0.1})}{N^2}+
\frac{(B_{1.0}-A_{0.1})\sqrt{\tilde y}}{\sqrt{-3bN\tilde x}}-\\
&-\frac{\sqrt{\tilde y}(BN_{1.0}-AN_{0.1})}{2\sqrt{-3b\tilde xN^3}}+
\frac{(B\tilde y_{1.0}-A\tilde y_{0.1})}{2\sqrt{-3bN\tilde x\tilde y}}
-\frac{\sqrt{\tilde y}(B\tilde x_{1.0}-A\tilde x_{0.1})}{2\sqrt{-3bN\tilde x^3}}=\\
&=\frac{50}{3}\Omega+ \frac{\sqrt{\tilde y}}{\sqrt{-3bN\tilde x}}\left(-\frac 65 N-
\frac{M}{2N}+\frac{5M}{2N} \right)=\frac{50}{3}\Omega+\frac{\sqrt{\tilde y}N}{\sqrt{-3bN\tilde x}}
\left(-\frac 65+2I_1\right)=0.
\endaligned
$$
The second condition true if  $\Omega=0$.

{\bf Theorem 3.} {\it An arbitrary equation (\ref{eq}) is equivalent to the Painleve III equation with three zero parameters  if and only if  the following conditions are true:
1) $F=0$ (\ref{F}), but $A\ne 0$ or $B\ne 0$ (\ref{alpha}),
2) $\Omega=0$ (\ref{Omega1}), (\ref{Omega2}) ,
3) $M\ne 0$ (\ref{M1}), (\ref{M2}), 4) $I_1=3/5.$
}

In this case, we can not write an explicit change of variables via
invariants of the equation because they are constants.

\section{Conclusion}
Thus we have found the explicite verification test for a second order
ODE to be equivalent to  Painleve equations I, II and III with three
zero parameters.
For the first two cases a point transformation of variables is found,
which is written using the differential invariants of the equation.

\section{Acknowlegments}

Author is grateful to Professor Rainer Picard (Institute of Analysis, TU Dresden) for the arrangement of conditions
for investigations and Professor Suleimanov Bulat Irekovich (Institute of Mathematics with Computing Centre  Ufa Science Centre Russian Academy of Science) for the information and useful discussions.

\section{Appendix}

Let us denote $K_{i.j}={\partial ^{i+j}K}/{\partial x^i\partial y^j}.$ 

The coordinates of the pseudovectorial field $\alpha$ are
$\alpha^1=B$, $\alpha^2=-A$,
where
\begin{equation}\label{alpha}
\aligned&\aligned A=P_{ 0.2}&-2Q_{ 1.1}+R_{ 2.0}+ 2PS_{ 1.0}+SP_{
1.0}-3PR_{ 0.1}-3RP_{ 0.1} -
3QR_{ 1.0} +6QQ_{ 0.1},\endaligned \\
\vspace{1ex}&\aligned B=S_{ 2.0}&-2R_{ 1.1}+Q_{ 0.2}- 2SP_{
0.1}-PS_{ 0.1}+3SQ_{ 1.0}+3QS_{ 1.0}+ 3RQ_{ 0.1}-6RR_{ 1.0}.
\endaligned\endaligned
\end{equation}

The first pseudoinvariant $F$ of the weight 5 is:
\begin{equation}\label{F}
3F^5=AG+BH,\qquad\text{where}
\end{equation}
$$ \aligned G&=-BB_{ 1.0}-3AB_{ 0.1}+4BA_{ 0.1}+
3SA^2-6RBA+3QB^2,\\
H&=-AA_{ 0.1}-3BA_{ 1.0}+4AB_{ 1.0}-
3PB^2+6QAB-3RA^2.
\endaligned
$$

The pseudoinvariant $N$ of the weight  2 in cases $A\ne 0$ and $B\ne 0$ respectively is: 
\begin{equation}\label{N}
 N =-\frac H{3A}, \qquad\qquad N =\frac G{3B}.
\end{equation}
The pseudoinvariant $M$ of the weight 4 in the case $A\ne 0$ is:
\begin{equation}\label{M1}
M=-\frac {12BN(BP+A_{ 1.0})}{5A}+BN_{ 1.0}+\frac {24}5BNQ+\frac
65NB_{ 1.0}+\frac 65NA_{ 0.1}-AN_{ 0.1}- \frac {12}5ANR
\end{equation}
and in the case $B\ne 0$ is:
\begin{equation}\label{M2}
M=-\frac {12AN(AS-B_{ 0.1})}{5B}-AN_{ 0.1}+\frac {24}5ANR- \frac
65NA_{ 0.1}-\frac 65NB_{ 1.0}+BN_{ 1.0}- \frac {12}5 BNQ.
\end{equation}
The pseudoinvariant  $\Omega$ of the weight 1 in the case $A\ne 0$ is:
\begin{equation}\label{Omega1}
\aligned \Omega &=\frac {2BA_{ 1.0}(BP+ A_{ 1.0})}{A^3}- \frac
{(2B_{ 1.0}+3BQ)A_{ 1.0}}{A^2}+\frac {(A_{ 0.1}-2B_{
1.0})BP}{A^2}-\\ &- \frac {BA_{ 2.0}+B^2 P_{ 1.0}}{A^2}+ \frac
{B_{ 2.0}}A+\frac {3B_{ 1.0}Q+3BQ_{ 1.0}- B_{ 0.1}P-BP_{
0.1}}{A}+Q_{ 0.1}- 2R_{ 1.0}
\endaligned
\end{equation}
and in the case  $B\ne 0$ is: 
\begin{equation}\label{Omega2}
\aligned \Omega &=\frac {2AB_{ 0.1}(AS- B_{ 0.1})}{B^3}- \frac
{(2A_{ 0.1}-3AR)B_{ 0.1}}{B^2}+\frac {(B_{ 1.0}-2A_{
0.1})AS}{B^2}+\\ &+ \frac {AB_{ 0.2}-A^2 S_{ 0.1}}{B^2}- \frac
{A_{ 0.2}}B+\frac {3A_{ 0.1}R+3AR_{ 0.1}- A_{ 1.0}S-AS_{
1.0}}{B}+R_{ 1.0}- 2Q_{ 0.1}.
\endaligned
\end{equation}

The pseudocovectorial field
 $\omega$ of the weight   -1  in the case $A\ne 0$ is:
$$
\aligned \omega_1=&\frac {12PR}{5A}-\frac {54}{25}\frac
{Q^2}A-\frac {P_{ 0.1}}A+ \frac {6Q_{ 1.0}}{5A}- \frac {PA_{
0.1}+BP_{ 1.0}+A_{ 2.0}}{5A^2}-\frac {2B_{ 1.0}P}{5A^2}+\\ &+
\frac {3QA_{ 1.0}-12PBQ}{25A^2}+
\frac {6B^2P^2+12BPA_{ 1.0}+6A_{ 1.0}^2}{25A^3},
\\
\omega_2=& \frac {-5BP_{ 0.1}+6BQ_{ 0.1}+12RBP}{5A^2} -\frac
{54}{25}\frac {BQ^2}{A^2}-\frac {2BB_{ 1.0}P+BA_{ 0.1}P+B^2P_{
1.0}+ BA_{ 2.0}}{5A^3}-\\ & -\frac {12B^2PQ}{25A^3}+ \frac {3BQA_{
1.0}}{25A^3}+
\frac {6BA_{ 1.0}^2+6B^3P^2+12B^2A_{ 1.0}P}{25A^4}
\endaligned
$$
and in the case $B\ne 0$ is:
$$
\aligned \omega_1=& \frac {5AS_{ 1.0}-6AR_{ 0.1}+12QAS}{5B^2}
-\frac {54}{25}\frac {AR^2}{B^2}+\frac {2AA_{ 0.1}S+AB_{
1.0}S+A^2S_{ 0.1}- AB_{ 0.2}}{5B^3}-\\ &-\frac
{12A^2SR}{25B^3}+\frac {3ARB_{ 0.1}}{25B^3}+
\frac {6AB_{ 0.1}^2+6A^3S^2-12A^2B_{ 0.1}S}{25B^4},\\
\omega_2=&\frac {12SQ}{5B}-\frac {54}{25}\frac {R^2}B+\frac {S_{
1.0}}B- \frac {6R_{ 0.1}}{5B}+ \frac {SB_{ 1.0}+AS_{ 0.1}-B_{
0.2}}{5B^2}+\frac {2A_{ 0.1}S}{5B^2}-\\ &- \frac {3RB_{
0.1}+12SAR}{25B^2}+ \frac {6A^2S^2-12B_{ 0.1}AS+6B_{
0.1}^2}{25B^3}.
\endaligned
$$

The pseudoinvariant  $\Theta$ of the weight -2 is given by:
\begin{equation}\label{Theta}
\Theta=\frac {\omega_1}A, \qquad \Theta=\frac {\omega_2}B.
\end{equation}

The pseudovectorial field $\theta$ of the weight -1 is:
\begin{equation}\label{theta}
\theta^1=\Theta_{ 0.1}-2\varphi_2\Theta, \qquad \theta^2=-\Theta_{
1.0}+2\varphi_1\Theta,
\end{equation}
where $\varphi_i$ in the case $A\ne 0$ are:
\begin{equation}\label{phi1}
\varphi_1=-3\frac {BP+A_{ 1.0}}{5A}+\frac 35Q, \;\;
\varphi_2=3B\frac {BP+A_{ 1.0}}{5A^2} -3\frac {B_{ 1.0}+A_{
0.1}+3BQ}{5A}+\frac 65R,
\end{equation}
and in the case $B\ne 0$   are:
\begin{equation}\label{phi2}
 \varphi_1=-3A\frac {AS-B_{ 0.1}}{5B^2}
-3\frac {A_{ 0.1}+B_{ 1.0}-3AR}{5B}-\frac 65Q,\quad
\varphi_2=3\frac {AS-B_{ 0.1}}{5B}-\frac 35R.
\end{equation}
The pseudoinvariant  L of the weight -4 is:
\begin{equation}\label{L}
\aligned
L=&
 \theta^1\theta^2(\theta^1_{
1.0}-\theta^2_{ 0.1})+ (\theta^2)^2\theta^1_{
0.1}-(\theta^1)^2\theta^2_{ 1.0}-\\
& -P(\theta^1)^3-3Q(\theta^1)^2\theta^2-3R\theta^1(\theta^2)^2-S(\theta^2)^3-\frac 12\Theta^2.
\endaligned
\end{equation}
The pseudoinvariant  $L_1$ of the weight -5 is:
\begin{equation}\label{L1}
L_1= L_{{1.0}}\theta^1+ L_{
0.1}\theta^2-4L(\varphi_1\theta^1+\varphi_2\theta^2).
\end{equation}
The pseudoinvariant $W$ of the weight -6 is:
\begin{equation}\label{W}
W=\nabla_{\theta}L_1= (L_1)_{
1.0}\theta^1+ (L_1)_{
0.1}\theta^2-5L_1(\varphi_1\theta^1+\varphi_2\theta^2).
\end{equation}
The pseudoinvariant  $V$ of the weight -3 is:
\begin{equation}\label{V}
V=\nabla_{\alpha}L_1= (L_1)_{
1.0}B- (L_1)_{
0.1}A-5L_1(B\varphi_1-A\varphi_2).
\end{equation}
Pseudovectorial field $\xi$ of the weight 3 is:
\begin{equation}\label{xi}
\xi
=-2\Omega\alpha-\gamma,
\end{equation}
in the case $A\ne 0$ field $\gamma$ is:
$$
\aligned
 \gamma^1=&-\frac {6BN(BP+A_{ 1.0})}{5A^2}+
\frac {18NBQ}{5A}+\\
&+\frac {6N(B_{ 1.0}+A_{ 0.1})}{5A} -N_{ 0.1}-\frac
{12}5NR-2\Omega B,
\endaligned
$$
$$
\gamma^2=-\frac {6N(BP+A_{ 1.0})}{5A}+N_{ 1.0}+\frac 65NQ+ 2\Omega
A,
$$
and in the case $B\ne 0$ is:
$$
\gamma^1=-\frac {6N(AN-B_{ 0.1})}{5B}-N_{ 0.1} +\frac 65NR-2\Omega
B,
$$
$$
\aligned
 \gamma^2=&-\frac {6AN(AS-B_{ 0.1})}{5B^2}+
\frac {18NAR}{5B}-\\
&-\frac {6N(A_{ 0.1}+B_{ 1.0})}{5B} +N_{ 1.0}-\frac
{12}5NQ+2\Omega A.
\endaligned
$$
The pseudoinvariant  $\Gamma$ of the weight 4 is:
\begin{equation}\label{Gamma}
\aligned \Gamma=&\frac {\gamma^1\gamma^2(\gamma^1_{
1.0}- \gamma^2_{ 0.1})}{M}+ \frac {(\gamma^2)^2\gamma^1_{ 0.1}-
(\gamma^1)^2\gamma^2_{ 1.0}}M+\\
&+\frac
{P(\gamma^1)^3+3Q(\gamma^1)^2\gamma^2+3R\gamma^1(\gamma^2)^2+
S(\gamma^2)^3}M.
\endaligned
\end{equation}


\bigskip
\begin{thebibliography}{1}

\bibitem{Liouville} R.~Liouville {\it Sur les invariants de certaines
equations differentielles et sur leurs applications} // J. de
L'Ecole Polytechnique. V.59. 1889. Pp.~7--76

\bibitem{Lie} S.Lie {\it Theorie der Transformationsgruppen III}
// Teubner Verlag. Leipzig. 1930.


\bibitem{Tresse1} A.~Tresse {\it Sur les invariants differenties des
groupes continus de transformations } // Acta Math. V 18. 1894. Pp.~
1--88. 

\bibitem{Tresse2} A.~Tresse {\it Determination des Invariants
ponctuels de l"Equation differentielle ordinaire de second ordre:
$y''=w(x,y,y')$} // Preisschriften der f\"rstlichen
Jablonowski'schen Gesellschaft XXXII. S.Hirzel.  Leipzig. 1896.


\bibitem{Cartan} E.~Cartan {\it Sur les varietes a connection
projective} // Bulletin de Soc. Math. de France, V. 52. 1924. Pp.~
205--241.




\bibitem{Thomsen} G.~Thomsen {\it {\"U}ber die topologischen
Invarianten der Differentialgleichung $y''=$
$f(x,y){y'}^3+g(x,y){y'}^2+h(x,y) y' +k(x,y) $} // Abhandlungen
aus dem mathematischen Seminar der Hamburgischen Universit\"at. V.
7. 1930. Pp. 301--328.

\bibitem{Grissom} C.Grissom, G.Thompson and G.Wilkens
{\it Linearisation of Second Order Ordinary Differential Equations via
Cartan's Equivalence Method} // Diff. Equations. V. 77. Pp. 1-15.
1989.

\bibitem{Rom} Yu.R.~Romanovskii {\it Calculation of local symmetries
of second order ordinary differential equations 
by means of Cartan's method of equivalence} // Manuskript, Pp. 1--20

\bibitem{BordagDruima} L.A.~Bordag and V.S.~Dryuma {\it Investigation of
dynamical systems using tools of the theory of invariants and
projective geometry} // NTZ-Preprnt 24/95 "addr Leipzig, 1995;
 Electronic archive at LANL (1997). solv-int \#9705006. Pp. 1--18.
 


\bibitem{Its} A.~R.~Its and V.~Yu.~Novokshenov {\it The
Isomonodromic Deformation Method in the Theory of Painleve
Equations} // Lecture Notes in Mathematics, Vol. 1191.
Springer-Verlag. New York/Berlin. 1986.

\bibitem{Sigur}  Mark J. Ablowitz, Harvey Segur {\it Solitons and the inverse scattering transform }// Mir, Moscow, 1987.


\bibitem{Gromak}	V.I.Gromak, N.A.Lukashevich {\it Analytical properties of solutions of Painleve equations} // Minsk. 1990. 

\bibitem{Sharipov1} Dmitrieva~V.~V., Sharipov~R.~A.{\it On the
point transformations for the second order differential equations}
// Electronic archive at LANL (1997). solv-int \#9703003. Pp.
1--14.


\bibitem{Sharipov2} Sharipov~R.~A. {\it On the point transformations
for the equation $y''=P+3\,Q\,y'+3\,R\,{y'}^2+S\,{y'}^3$}
// Electronic archive at LANL (1997). solv-int \#9706003. Pp.
1--35.


\bibitem{Sharipov3} Sharipov~R.~A.{\it Effective procedure of point
classification for the equations $y''=P+3\,Q\,y'+3\,R\,{y'}^2
+S\,{y'}^3$} // Electronic archive at LANL (1998). Math\.DG
\#9802027. Pp. 1--35.






\bibitem{Lamb} N.Kamran, K.G.Lamb \& W.F.Shadwick {\it The local
equivalence problem for $d^2y/dx^2=F(x,y,dy/dx)$ and the Painleve
transcendents} // J.Differential geometry. V. 22. 1985.  Pp. 139-150.

 \bibitem{Sokolov} A.V.Bocharov, V.V.Sokolov, and S.I.Svinolupov 
 {\it On some equivalence problem for diffrential equations} // Preprint ESI-54, International Erwin Schr\"odinger Institute for Mathematical Physics, Wien, Austria, 1993, pp.1-12.
 
 \bibitem{Babich} M.V.~Babich and L.A.~Bordag {\it Projective
Differential Geometrical Structure ot the Painleve Equations} //
J. of Diff.Equations. V. 157 (2).   1999. Pp. 452--485.
 
 \bibitem{Hietarinta} Hietarinta J. and Dryuma V. {\it Is my ODE is Painleve equation in disguise? } // J.of Nonlin. Math. Phys. 2002. 9(1). pp.~67-74. 
 
\bibitem{Dridi} Raouf Dridi {\it On the geometry of the first and second Painleve equations} //
J. Phys. A: Math. Theor. V. 42. 2009. pp.1-9.



\bibitem{Kartak} Kartak V.V. {\it Solution of the equivalence problem for Painleve I and II equations} // VINITI No 612-B 2006.


\bibitem{Kamke} E.Kamke {\it Handbook of ordinary differential equations} //  Moscow. Nauka. 1976.





\end{thebibliography}
\end{document}